\documentclass[10pt]{article}
\usepackage{mathrsfs}
\usepackage{stmaryrd}
\usepackage{amsfonts,amsmath,amssymb,amscd}
\usepackage{shadow}
\usepackage{indentfirst,graphics,epsfig,psfrag}
\usepackage{color}
\usepackage{longtable}
\usepackage{pstricks,multido}
\usepackage{hyperref}
\usepackage{qtree}
\input{epsf}
\usepackage{ifpdf}
\usepackage{enumerate}
\usepackage{appendix}
\usepackage{enumerate}
\usepackage{lineno}
\allowdisplaybreaks

\newtheorem{thm}{Theorem}[section]

\newtheorem{prop}[thm]{Proposition}

\def\qed{\hfill \rule{4pt}{7pt}}

\def\pf{\noindent {\it{Proof.} \hskip 2pt}}

\numberwithin{equation}{section}

\begin{document}
\begin{center}
{\large\bf
{The minimum size of graphs with given rainbow index}
}
 \end{center}

\begin{center}
{Thomas Y.H. Liu}
\vskip 2mm
Department of Foundation Courses\\
Southwest Jiaotong University\\ Emeishan,
Sichuan 614202, P.R. China\\

\vskip 2mm
E-mail\,$:$thomasliuyh@sina.com
\end{center}

\begin{abstract}
The concept of $k$-rainbow index $rx_k(G)$ of a connected graph $G$,
introduced by Chartrand, Okamoto and Zhang, is a natural
generalization of the rainbow connection number. Let $t(n,k,\ell)$
denote the minimum size of a connected graph $G$ of order $n$ with
$rx_k(G)\leq \ell$, where $2\leq \ell\leq n-1$ and $2\leq k\leq n$.
In this paper, we obtain some exact
values and some upper bounds for $t(n,k,\ell)$.
\end{abstract}

\noindent {\bf Keywords}: edge coloring, rainbow connection, rainbow index, rainbow $S$-tree

\noindent {\bf AMS  Subject Classifications}: 05C15, 05C35, 05C40

\section{Introduction}

All graphs considered in this paper are finite, undirected and
simple. We follow the notation and terminology of Bondy and Murty
\cite{Bondy}, unless otherwise stated.

An edge-colored graph $G$ is \emph{rainbow connected} if every two
vertices are connected by a path satisfying no two edges on the path have the same color.  The
minimum number of colors required to make a graph $G$ rainbow connected is
called \emph{the rainbow connection number}, denoted $rc(G)$.   The notion of rainbow
connection in graphs was introduced by Chartrand et
al. \cite{Chartrand}.

For an edge-colored nontrivial connected graph $G$ of order $n$, A tree $T$ in $G$ is called a \emph{rainbow tree} if no two edges of $T$ are colored
the same. For $S\subseteq V(G)$ and $|S|\geq2$, a \emph{rainbow $S$-tree} is a
rainbow tree $T$ such that $S\subseteq V(T)$. For a fixed integer $k$ with $2\leq
k\leq n$, an edge-coloring $c$ of $G$ is called a \emph{$k$-rainbow
coloring} if for every $k$-subset $S$ of $V(G)$ there exists a
rainbow $S$-tree. The minimum number of colors that are needed in
a $k$-rainbow coloring of $G$ is called the \emph{$k$-rainbow index}
of $G$, denoted $rx_k(G)$. Clearly, when $k=2$, $rx_2(G)$ is the rainbow connection number $rc(G)$ of $G$.
Note that $k$-rainbow index, defined by Chartrand et al.\cite{Chartrand2}, is a generalization of rainbow connection number. The study about rainbow connection has been extensively researched, we refer to \cite{clrty,lll,lms,LiSun,LiSun1} for example.

In this paper, motivated by a recent paper of
Schiermeyer \cite{Schiermeyer} on the minimum size of rainbow $k$-connected graphs, where a graph $G$ is called \emph{rainbow k-connected} if there is an
edge colouring of $G$ with $k$ colours such that $G$ is rainbow connected, we study the minimum size of a graph $G$ such that $G$ has a $k$-rainbow coloring using a fixed number of colors. To be more specific, let $t(n,k,\ell)$
be the minimum size of a connected graph $G$ of order $n$ with
$rx_k(G)\leq \ell$, where $2\leq \ell\leq n-1$ and $2\leq k\leq n$. Observe that
\[
t(n,k,1)\geq t(n,k,2)\geq \ldots \geq t(n,k,n-1).
\]
Our main objective is to give some exact values and some upper bounds for $t(n,k,\ell)$ when $k$ and $\ell$ take specific values.

\section{Main Results}\label{mainresults}
In this section, we mainly concern about some  exact values and some upper bounds for $t(n,k,\ell)$, when $k=3$.

\begin{prop}\label{lmthm2-1} Let $n\geq3$ be a positive integer. Then

$(1)$
\[ t(n,3,2)=
\begin{cases}
\,\,2,&\mathrm{if}~n=3;\\
\,\,4,&\mathrm{if}~n=4;\\[1mm]
{{5}\choose{2}},&\mathrm{if}~n=5.
\end{cases}
\]
Furthermore, when $n\geq6$, there does not
exist a connected graph $G$
such that $rx_3(G)\leq2$.

$(2)$ \[
  t(n,3,3)=
    \begin{cases}
             \,\,2,&\mathrm{if}~n=3; \\
             \,\,3,&\mathrm{if}~n=4;\\
             \,\,5,&\mathrm{if}~n=5.
   \end{cases}
  \]

Furthermore, when $n\geq6$,
\[
  t(n,3,3)\leq
    \begin{cases}
             \,\,\frac{n^2}{4},&\text{if}~n \text{ is even}; \bigskip\\
             \,\,\frac{(n+3)(n-1)}{4},&\text{if}~n \text{ is odd}.
   \end{cases}
  \]
  \end{prop}
To prove Proposition \ref{lmthm2-1}, we need the following results.

\begin{thm}[\mdseries{\cite[Proposition 1.3]{Chartrand2}}]\label{czthm1}
Let $T$ be a tree of order $n\geq 3$. For each integer $k$ with $3\leq k\leq n$, $$rx_k(T)=n-1.$$
\end{thm}

\begin{thm}[\mdseries{\cite[Theorem 2.1]{Chartrand2}}]\label{cztheorem1}
For integers $k$ and $n$ with $3\leq k\leq n$,
\begin{equation*}
 rx_k(C_n)=
   \begin{cases}
      n-2, & \text{if $k = 3$ and $n\geq4$}; \\
     n-1, & \text{if $k = n = 3$ or $4\leq k\leq n$.}
    \end{cases}
\end{equation*}
\end{thm}

\begin{thm}[\mdseries{\cite[Theorem 2.2 and Theorem 2.3]{Chartrand2}}]\label{ctthm2}
If $G$ is a unicyclic graph of order $n\geq 3$ and girth $g\geq 3$ that is not a cycle, then
\begin{equation*}
 rx_3(G)=
   \begin{cases}
      n-2, & \text{$g\geq4$}; \\
     n-1, & \text{$g=3$};
    \end{cases}
\end{equation*}
where a unicyclic graph means a connected graph containing exactly one cycle.
\end{thm}

\begin{thm}[\mdseries{\cite[Theorem 3 ]{CLYZ}}]\label{clyzthm1}
Let $G$ be a connected graph of order $n$. Then $rx_3(G)=2$ if and
only if $G=K_5$ or $G$ is a 2-connected graph of order 4 or $G$ is
of order 3.
\end{thm}
\begin{thm}[\mdseries{\cite[Theorem 5]{CLYZ}}]\label{clyz5} For each integer $r$ with $r\geq3$, $rx_3(K_{r;r})=3$.
\end{thm}

Using  Theorem \ref{czthm1}, Theorem \ref{cztheorem1}, Theorem \ref{ctthm2}, Theorem \ref{clyzthm1} and Theorem \ref{clyz5}, it is easy to prove Proposition \ref{lmthm2-1}, we sketch it as follows.

\noindent
\textit{Proof of Proposition \ref{lmthm2-1}.}
$(1)$ The result follows by Theorem \ref{clyzthm1}.

$(2)$ It is clearly true when $n=3$ or $n=4$.
Let $G$ be a graph with $5$ vertices and $rx_3(G)\leq3$. If $G$ is a tree, then from  Theorem \ref{czthm1}, $rx_3(G)=4$, a contradiction. Thus $G$ must contain a cycle, from Theorem \ref{cztheorem1} and Theorem \ref{ctthm2}, it is easy to deduce that the minimum size of $G$ is $5$.

When $n\geq6$, if $n$ is even, then from Theorem \ref{clyz5}, we have $rx_3(K_{\frac{n}{2};\frac{n}{2}})=3.$
Therefore, it is clear that
$t(n,3,3)\leq \frac{\,\,n^2}{4}$.

If $n$ is odd, let
$U=\{u_1,u_2,\ldots,u_\frac{n-1}{2}\}$
, $V=\{v\}$
and $W=\{w_1,w_2,\ldots,w_{\frac{n-1}{2}}\}$. Then the join $G[U,W]\vee v$ of a complete bipartite graph $G[U,W]$ and the isolate vertex $v$ is denoted by $H$ with $V(H)=U\cup V\cup W$ and $E(H)=\{(x,v)\,|\,x\in U\cup W \}\cup\{(x,y)\,|\,x\in U \text{ and } y\in W\}$.
Now let us prove that $$rx_3(H)=3.$$
Similar to the proof of Theorem \ref{clyz5}, define a coloring $c$:
$E(H)\rightarrow \{1,2,3\}$ as follows:
 \begin{equation}
c(u_iw_j)=
   \begin{cases}
      1, & 1\leq i=j\leq  \frac{n-1}{2}; \\
     2, & 1\leq i< j\leq  \frac{n-1}{2}; \\
    3, & 1\leq j< i\leq   \frac{n-1}{2};
    \end{cases}
\end{equation}
and
\begin{equation}
c(xv)=
   \begin{cases}
     2, & \text{if $x\in U$}; \\
    3, & \text{if $x\in W$}.
    \end{cases}
\end{equation}
 Now, we show that
 $c$ is a $3$-rainbow coloring of $K_{H}$. Let $S\subseteq V(K_{H})$ with $|S|=3$.
Then we consider the following two cases.

{\bf Case 1.}~~$S\subseteq U \cup W$.

Under this case, since $G[U,W]$ is $K_{\frac{n-1}{2};\frac{n-1}{2}}$ and
$rx_3({K_{\frac{n-1}{2};\frac{n-1}{2}}})=3$, then
it is routine to verify that there is a rainbow $S$-tree,
see \cite{CLYZ}.

{\bf Case 2.}~~$v\in S$.

Let $S=\{x,y,v\}$. Then if $x\in U$ and $y\in W$,
$T=\{vx,vy\}$ is a rainbow
$S$-tree. If $x,y\in U$ with $x=u_i$ and $y=u_j$, then $T=\{vu_i,vw_j,u_jw_j\}$ is a rainbow
$S$-tree. If $x,y\in W$ with $x=w_i$ and $y=w_j$, then $T=\{vw_i,vu_j,u_jw_j\}$ is a rainbow
$S$-tree.

Therefore $$rx_3(H)\leq3.$$
Moreover, it is clearly that $rx_3(H)\geq3.$
Thus$$rx_3(H)=3.$$
Hence $$t(n,3,3)\leq \frac{(n-1)^2}{4}+(n-1)=\frac{(n+3)(n-1)}{4}.$$

\qed

\begin{thm}\label{lmthm2-2}
Let $n\geq3$ be a positive integer. Then
$$t(n,3,4)\leq {n \choose 2}-n+1.$$
\end{thm}
\pf
Let $G$ be a graph such that $\overline{G}$ is a union of a cycle of
order $n-1$ and an isolated vertex. Let $w$ be the isolated
vertex in $\overline{G}$. Then $d_G(w)=n-1$. Set $$V(G)\setminus
w=\{v_1,v_2,\ldots,v_{n-1}\}.$$
Let $n-1=3r+t$ where $0\leq t\leq 2$. For $0\leq j\leq r, 1\leq i\leq n-1$, set
$$ X_1=\{v_i\,|\,i=3j+1\},$$
$$X_2=\{v_i\,|\,i=3j+2\}\,$$
and
$$X_3=\{v_i\,|\,i=3j+3\}.$$
Then $G[X_1]$ is a clique or a graph
obtained from a clique of order $|G[X_1]|$ by deleting one edge,
both $G[X_2]$ and $G[X_3]$ are cliques.

In order to  show $rx_3(G)\leq 4$, we provide an edge-coloring $c:
E(G)\longrightarrow \{1,2,3,4\}$ defined by
$$
c(e)=\left\{
\begin{array}{ccc}
1,&\mathrm{if}~e\in E_G[w,X_1]&\cup~~ E(G[X_3]);\qquad\qquad\qquad\qquad\qquad\qquad\\[1mm]
2,&\mathrm{if}~e\in E_G[w,X_2]&\cup~~ E(G[X_1]);\qquad\qquad\qquad\qquad\qquad\qquad\\[1mm]
3,&\mathrm{if}~e\in E_G[w,X_3]&\cup~~ E(G[X_2]);\qquad\qquad\qquad\qquad\qquad\qquad\\[1mm]
4,&~~\mathrm{if}~e\in E_G[X_1,X_2]&\cup~~E_G[X_1,X_3]~\cup~ E_G[X_2,X_3].\,\qquad\,\,\qquad
\end{array}
\right.
$$
Clearly, $c(wv_{n-1})=1$ if and only if $n-1=3r+1$ for some positive integer $r$. It suffices
to show that there exists a rainbow $S$-tree for any $S\subseteq
V(G)$ with $|S|=3$.

\textbf{Case 1}: $w\in S$.

Without loss of generality, let $$S=\{w,v_i,v_j\} \ (i<j).$$ If
$c(wv_i)\neq c(wv_j)$, then the tree induced by the edge set
$\{wv_i,wv_j\}$ is a rainbow $S$-tree, as desired. So we assume
$$c(wv_i)=c(wv_j).$$ If $$c(wv_i)=c(wv_j)=1,$$ then the tree induced by the edge set
$$\{wv_{i+2},v_{i}v_{i+2},wv_j\}$$ is a rainbow $S$-tree with colors
$\{1,3,4\}$. When
$$c(wv_i)=c(wv_j)=2$$
or
$$c(wv_i)=c(wv_j)=3,$$
it is easy to verify that the tree induced by the edge set
$$\{wv_{i+2},v_{i}v_{i+2},wv_j\}$$ is a rainbow $S$-tree, as desired.

\textbf{Case 2}: $w\notin S$.

Without loss of generality, let $S=\{v_i,v_j,v_k\} \ (i<j<k)$. Firstly, if
$$c(wv_i)\neq c(wv_j)\neq c(wv_k),$$ then clearly the tree induced by the edge set
$$\{wv_i, wv_j, wv_k\}$$
is the $S$-tree with colors $\{1,2,3\}$.

Secondly, we consider the case $$c(wv_i)=c(wv_j)=c(wv_k).$$ Suppose
$$c(wv_i)=c(wv_j)=c(wv_k)=1,$$ then $$v_i,v_j,v_k\in X_1\text{ and } v_j \text{ is adjacent to } v_k.$$
Since $$c(v_iv_{i+2})=4,  c(v_jv_k)=2, c(wv_j)=1$$ and
$$c(wv_{i+2})=3.$$
Therefore, the tree induced by the edge set
$$\{wv_{i+2},v_{i}v_{i+2},v_jv_k,wv_j\}$$ is a rainbow $S$-tree with
colors $\{1,2,3,4\}$, as desired. Similarly, when $$c(wv_i)=c(wv_j)=c(wv_k)=2$$
or
$$c(wv_i)=c(wv_j)=c(wv_k)=3.$$
It is routine to verify that the tree induced by the edge set
$$\{wv_{i+2},v_{i}v_{i+2},v_jv_k,wv_j\}$$ is a rainbow $S$-tree with
colors $\{1,2,3,4\}$.

Finally, we consider the case that only two edges in $\{wv_i,wv_j,wv_k\}$
receive the same color under the coloring $c$. This implies that there exist $$v_l,v_m\in S=\{v_i,v_j,v_k\}$$ satisfying $v_l,v_m\in X_p$ for some $p\in\{1,2,3\}$ and the left element $v_t\in S\setminus\{v_l,v_m\}$ is in $X_q$ with $q\neq p$. Moreover, $v_t$ must be adjacent to $v_l$ or $v_m$, without loss of generality, set $v_t$ adjacent to $v_m$.
Let $$x\in\{1,2,3\}\setminus\{p,q\} \text{ and } v_y\in X_x.$$
Then $$c(v_tv_m)=4,c(wv_t)=q,c(wv_y)=x \text{ and } c(wv_l)=p.$$
Thus the tree induced by the edge set
$$\{v_tv_m,wv_t,wv_y,wv_l\}$$
is a rainbow $S$-tree with colors \{1,2,3,4\}.

From the above arguments, we conclude that $rx_3(G)\leq 4$ and
$t(n,3,4)\leq {n \choose 2}-n+1$.

\qed

The join $C_n\vee K_1$ of a cycle $C_n$ and a single vertex is referred to as a \emph{wheel} with $n$ \emph{spokes}, denoted $W_n$, see \cite{Diestel}.
The $3$-rainbow index of $W_n$ was given by Chen, Li, Yang and Zhao as follows.

\begin{thm}[\mdseries{\cite[Theorem 7]{CLYZ}}]\label{clyz3}
For $n\geq 3$, the $3$-rainbow index of the wheel $W_n$ is
\begin{equation}
 rx_3(W_n)=
   \begin{cases}
     2, & n=3; \\
     3, & 4\leq n\leq 6;\\
     4, & 7\leq n\leq 16;\\
     5, & n\geq 17.
    \end{cases}
\end{equation}
\end{thm}
Thus
 $$rx_3(W_n)\leq5.$$
Since $|E(W_n)|=2n$, the following result is true.
\begin{prop}\label{lmthm2-3} Let
$n\geq3$ be an integer. Then
$$t(n,3,5)\leq 2n-2.$$
\end{prop}

\begin{thm}\label{lmthm2-4}Let $n\geq3$ be an integer. Then
$$t(n,3,6)\leq 2n-3.$$
\end{thm}
\pf
When $3\leq n\leq 4$, it is true clearly. When $n\geq 5$, let $u$ be an isolated vertex and $v$ be the center of $W_{n-2}$. $G$ is the graph by adding an edge between $u$ and $v$, that is,
$$V(G)=V(W_{n-2})\cup\{u\}$$
and
$$E(G)=E(W_{n-2})\cup\{uv\}.$$
Since from Theorem \ref{clyz3} $$rx_3(W_{n-2})\leq 5,$$
it is easy to know that
$$rx_3(G)\leq 6.$$
Therefore, for any positive integer $n\geq 3$, we have
$$t(n,3,6)\leq 2n-3,$$
since $|E(G)|=2n-3.$
\qed

\begin{thm}\label{lmthm2-5}
Let $n$ and $\ell$ be positive integers satisfying $7\leq \ell\leq \frac{n-1}{2}$. Then
$$
t(n,3,\ell)\leq n+t+{t\choose 2}\Big (n-2-\Big
\lfloor\frac{n-2}{\ell-3}\Big \rfloor(\ell-3)\Big)+{t+\lfloor\frac{n-2}{\ell-3}\rfloor-\big
\lceil\frac{n-2}{\ell-3}\big \rceil\choose
2}\Big (\ell+1-n+\Big \lfloor\frac{n-2}{\ell-3}\Big
\rfloor(\ell-3)\Big),
$$
where $t=\big \lceil\frac{n-2}{\ell-3}\big \rceil$.
\end{thm}

\pf
Let $$t=\big
\lceil\frac{n-2}{\ell-3}\big \rceil.$$
For $1\leq i\leq t-1$, let
$$Q_i=uv_{i,1}v_{i,2}\ldots v_{i,\ell-3}w$$
be a path of length $\ell-2$.
If $\big
\lceil\frac{n-2}{\ell-3}\big \rceil=\lfloor\frac{n-2}{\ell-3}\rfloor$, then let
$$Q_t=uv_{t,1}v_{t,2}\ldots v_{t,\ell-3}w$$
be a path of length $\ell-2$.
Otherwise let
$$Q_t=uv_{t,1}v_{t,2}\ldots v_{t,s}w$$
be a path of length $s+1$, where
$s=n-\lfloor\frac{n-2}{\ell-3}\rfloor(\ell-3)-2$.

Let $H$ be the union of $Q_1,Q_2,\ldots,Q_t$.
Let $G$ be the graph obtained from $H$ by adding the edges in
$$
\big \{uw\big\}\cup \big\{v_{i_1,j}v_{i_2,j}\,|\,1\leq i_1\neq
i_2\leq t, 1\leq j \leq s \big \}\cup \big
\{v_{i_1,j}v_{i_2,j}\,|\,1\leq i_1\neq i_2\leq t, s+1\leq j\leq
\ell-3 \big \}.
$$
Note that for each $j \ (1\leq j\leq s)$ the graph induced by the
vertex set $$\{v_{i,j}\,|\,1\leq i\leq t\}$$ is a complete graph of
order $t$. For each $j \ (s+1\leq j\leq \ell-3)$, the graph induced
by the vertex set $$\{v_{i,j}\,|\,1\leq i\leq t-1\}$$ is a complete
graph of order $t-1$.

In order to prove $rx_3(G)\leq \ell$, we provide  an edge-coloring $c:
E(G)\longrightarrow \{1,2,\ldots,\ell\}$ as follows:
If $\big
\lceil\frac{n-2}{\ell-3}\big \rceil=\lfloor\frac{n-2}{\ell-3}\rfloor$, $c(e)$ is defined to be
\begin{equation*}
 c(e)=
   \begin{cases}
     1, &\text{if $e=uv_{i,1}$ for $1\leq i\leq t$};\\
     j, & \text{if $e=v_{i,j-1}v_{i,j}$ for  $1\leq i\leq t$ and $2\leq j\leq \ell-3$}; \\
     \ell-2,& \text{if $e=v_{i,\ell-3}w$ for $1\leq i\leq t$};\\
     \ell-1, &\text{if $e=uw$};\\
     \ell,&\text{others}.
    \end{cases}
\end{equation*}
Otherwise,  $c(e)$ is defined to be

\begin{equation*}
 c(e)=
   \begin{cases}
     1, &\text{if $e=uv_{i,1}$, for $1\leq i\leq t$};\\
     j, & \text{if $e=v_{i,j-1}v_{i,j}$ or $e=v_{t,p-1}v_{i,q}$, for  $1\leq i\leq t-1$, $2\leq j\leq \ell-3$ and $2\leq q\leq s$}; \\
     \ell-2,& \text{if $e=v_{i,\ell-3}w$ or $e=v_{t,s}w$, for $1\leq i\leq t-1$};\\
     \ell-1, &\text{if $e=uw$};\\
     \ell,&\text{others}.
    \end{cases}
\end{equation*}
To show $rx_3(G)\leq \ell$, it suffices to prove that there exists a
rainbow $S$-tree for any $S\subseteq V(G)$ with $|S|=3$.
Set $$S=\{x,y,z\}.$$

\textbf{Case 1.} $S\subset V(Q_i)$ for some $1\leq i\leq t.$

It is easy to know that $Q_i$ is the required rainbow $S$-tree.

\textbf{Case 2.} Two of vertices in $S$ are on $Q_i$ for some $1\leq i\leq t$.

Without loss of generality, let $x,y\in V(Q_i)$ and $z\in V(Q_j)$, where $j\neq i$ and $z$ is neither $u$ nor $w$.
Suppose $$x=v_{i,p},y=v_{i,q}\text{ and }z=v_{j,r},$$
where $p\leq q$.
If $p\leq q \leq r$, then the path $$P=v_{i,p}v_{i,p+1}\ldots v_{i,q-1}v_{i,q}v_{j,q}v_{j,q+1}\ldots v_{j,r-1}v_{j,r}$$
is a rainbow $S$-tree. If $r\leq p \leq q$ or $p\leq r \leq q$, it is similar, so we omit.

\textbf{Case 3.} For any $1\leq i\leq t$, $S\cap V(Q_i)=1$.

Suppose
$$x=v_{i,p},y=v_{j,q} \text{ and } z=v_{k,r},$$
where $i,j,k,p,q$ and $r$ are positive integers satisfying $1\leq i\neq j\neq k\leq t$ and $p\leq q\leq r$. If $$\big
\lceil\frac{n-2}{\ell-3}\big \rceil=\lfloor\frac{n-2}{\ell-3}\rfloor\text{ or }k\neq t,$$ then the path
$$P'=v_{k,r}v_{k,r+1}\ldots v_{k,\ell-3}wuv_{i,1}v_{i,2}\ldots v_{i,p}v_{j,p}v_{j,p+1}\ldots v_{j,q}$$
is a rainbow $S$-tree.
If $$\big
\lceil\frac{n-2}{\ell-3}\big \rceil\neq\lfloor\frac{n-2}{\ell-3}\rfloor\text{ and }k=t,$$ then the path
$$P'=v_{k,r}v_{k,r+1}\ldots v_{k,s}wuv_{i,1}v_{i,2}\ldots v_{i,p}v_{j,p}v_{j,p+1}\ldots v_{j,q}$$
is a rainbow $S$-tree.

Therefore, $$rx_3(G)\leq \ell.$$
Since
$$
|E(G)|=n+t+{t\choose 2}\Big (n-2-\Big
\lfloor\frac{n-2}{\ell-3}\Big \rfloor(\ell-3)\Big)+{t+\lfloor\frac{n-2}{\ell-3}\rfloor-\big
\lceil\frac{n-2}{\ell-3}\big \rceil\choose
2}\Big (\ell+1-n+\Big \lfloor\frac{n-2}{\ell-3}\Big
\rfloor(\ell-3)\Big),
$$
the theorem is confirmed.
\qed

A \emph{rose graph} $R_{p}$ with $p$ \emph{petals} (or $p$-\emph{rose
graph}) is a graph obtained by taking $p$ cycles with just a vertex
in common. The common vertex is called the \emph{center} of $R_{p}$.
If the length of each cycle is exactly $q$, then this rose graph
with $p$ petals is called a $(p,q)$-\emph{rose graph}, denoted $R_{p,q}$.

\begin{thm}\label{lmthm2-6}
For $\frac{n}{2}\leq \ell\leq n-3$, $$t(n,3,\ell)\leq 2n-\ell-1.$$
\end{thm}
\pf
Let $G$ be a graph obtained from a $(n-\ell,3)$-rose graph
$R_{n-\ell,3}$ and a path $P_{2\ell-n}$ by identifying the center of
the rose graph and one endpoint of the path. Clearly,
$$|V(G)|=(2(n-\ell)+1)+(2\ell-n-1)=n$$ and $$|E(G)|=3(n-\ell)+(2\ell-n-1)=2n-\ell-1.$$

Let $w_0$ be the center of $R_{n-\ell,3}$, and let $C_i=w_0v_iu_iw_0
\ (1\leq i\leq n-\ell)$ be a cycle of $R_{n-\ell,3}$. Let
$P_{2\ell-n}=w_0w_1\ldots w_{2\ell-n-1}$ be the path of order
$2\ell-n$. To show that $rx_3(G)\leq \ell$, we provide an
edge-coloring $c: E(G)\longrightarrow \{1,2,\ldots,\ell\}$ defined
by
$$
c(e)=\left\{
\begin{array}{ll}
i,&\mathrm{if}~e=w_0u_i~\mathrm{or}~e=w_0v_i,\text{ for } 1\leq i\leq n-\ell;\\
n-\ell+i,&\mathrm{if}~e=w_{i-1}w_i,\text{ for } 1\leq i\leq 2\ell-n-1;\\
\ell, &\mathrm{if}~e=u_iv_i,\text{ for } 1\leq i\leq n-\ell;
\end{array}
\right.
$$
It suffices to show that there exists a rainbow $S$-tree for any
$S\subseteq V(G)$ and $|S|=3$. Set
$$S=\{x,y,z\}.$$

\textbf{Case 1} $w_0\in S$.

Set $x=w_0$. If $$y,z\in V(C_i)\setminus\{w_0\}$$ for some $1\leq i \leq n-\ell,$ then the graph induced by the edge set $\{xy,yz\}$ is a rainbow $S$-tree. If $$y\in V(C_i)\setminus\{w_0\},z\in V(C_j)\setminus\{w_0\}$$
for some $1\leq i\neq j\leq n-\ell,$ then the graph induced by the edge set $\{xy,xz\}$ is a rainbow $S$-tree.
If $$y\in V(C_i)\setminus\{w_0\},z\in V(P_{2\ell-n})\setminus\{w_0\}$$ for some $1\leq i\leq n-\ell,$ then the path $yxw_1w_2\ldots z$ is a rainbow $S$-tree. If $$y,z\in V(P_{2\ell-n})\setminus\{w_0\},$$
then $P_{2\ell-n}$ is a rainbow $S$-tree.

\textbf{Case 2} $S\subset V(P_{2\ell-n})\setminus\{w_0\}$.

Clearly, $P_{2\ell-n}$ is a rainbow $S$-tree.

\textbf{Case 3} $S\subset V(R_{n-\ell,3})\setminus\{w_0\}$.

If $$x,y\in V(C_i) \text{ and } z\in V(C_j)$$
for some $1\leq i\neq j\leq n-\ell$,
then the path $zw_0xy$ is a rainbow $S$-tree.
If $$x\in V(C_i), y\in V(C_j) \text{ and } z\in V(C_k)$$
for some $1\leq i\neq j\neq k\leq n-\ell$,
then the star induced by the edge set $\{w_0x,w_0y,w_0z\}$ is a rainbow $S$-tree.

\textbf{Case 4} $$|S\cap (V(R_{n-\ell,3})\setminus\{w_0\})|+|S\cap (V(P_{2\ell-n,3})\setminus\{w_0\})|=3$$ and $$|S\cap (V(R_{n-\ell,3})\setminus\{w_0\})|\cdot|S\cap (V(P_{2\ell-n,3})\setminus\{w_0\})|\neq0.$$
Similarly, it is routine to verify that there must exist a rainbow $S$-tree in $G$.

From the above arguments, we achieve that $$rx_3(G)\leq \ell$$ and
$$t(n,3,\ell)\leq 2n-\ell-1, \text{ for } \frac{n}{2}\leq \ell\leq n-3.$$
\qed

\begin{prop}\label{lmthm2-7}
Let $n\geq4$ be a positive integer. Then

$(1)$ $t(n,3,n-2)=n$;

$(2)$ $t(n,3,n-1)=n-1$.
\end{prop}
To prove Proposition \ref{lmthm2-7}, we need the following results.
\begin{thm}[\mdseries{\cite[Theorem 3]{LSYZ}}]\label{thm4}
Let $G$ be a connected graph of order $n$. Then $rx_3(G)=n-1$ if and only if $G$ is a tree or $G$ is a unicyclic graph with girth 3.
\end{thm}
Combining Theorem \ref{cztheorem1} and Theorem \ref{thm4}, Proposition \ref{lmthm2-7} is clearly true.

We conclude this paper with the following theorem about an upper bound of $t(n,n-1,n-2)$.
\begin{thm}\label{1mthm2-8}
$t(n,n-1,n-2)\leq 2n-4$.
\end{thm}
\pf
Let $G=K_{2,n-2}=G[X,Y]$, where set $X=\{u,w\}$ and
$Y=\{v_1,v_2,\ldots,v_{n-2}\}$. Then we give  an edge coloring $c$ of $G$ as follows:
\begin{equation*}
 c(e)=
   \begin{cases}
     i, & \text{$e=uv_i$, for $1\leq i\leq n-2$}; \\
     n-1-i, & \text{$e=wv_i$, for $1\leq i\leq n-2$ }.
    \end{cases}
\end{equation*}
We prove that $G$ is rainbow
$3$-tree-connected under this coloring. Let $$S\subseteq V(G)\text{ and } |S|=n-1.$$
If $S=V(G)\setminus \{u\}$, then the tree $T$ induced by the
edge set $$\{wv_1,wv_2,\ldots,wv_{n-2}\}$$ is a rainbow $S$-tree.
If $S=V(G)\setminus \{w\}$, there exists a rainbow $S$-tree similarly. If
$S=V(G)\setminus \{v_{i}\}$, then the tree $T$ induced by
the edge set
$$\{uv_1,uv_2,\ldots,uv_{i-1},uv_{i+1},\ldots,uv_{n-2},wv_{n-i-1}\}$$
is a rainbow $S$-tree. Therefore we conclude
$$rx_{n-1}(G)\leq n-2.$$ Since $e(G)=2n-4$, it
follows that $t(n,n-1,n-2)\leq 2n-4$. \qed

\vspace{0.5cm}
 \noindent{\bf Acknowledgments.}  I am indebted to the referee for a careful reading of the paper
and many useful suggestions. I am also very grateful for the valuable comments and discussion with Y. Mao and Y. Zhao.

\end{document}